\documentclass[10pt]{article}

\usepackage{amssymb}
\usepackage{amsthm}
\usepackage{amsmath}

\newtheorem{theorem}{Theorem}
\newtheorem{lemma}[theorem]{Lemma}

\theoremstyle{definition}

\theoremstyle{remark}

\numberwithin{equation}{section}
\numberwithin{theorem}{section}

\def \tilde{\widetilde}

\def \p{\partial}

\newcommand{\st}[1]{\ensuremath{^{\scriptstyle \textrm{#1}}}}

\newcommand{\CC}{{\mathbb C}}
\newcommand{\NN}{{\mathbb N}}

\newcommand{\ZZ}{{\mathbb Z}}

\newcommand{\D}{\mathcal{D}}

\renewcommand{\O}{\mathcal{O}}

\newcommand{\fg}{{\mathfrak g}}

\newcommand{\fo}{{\mathfrak o}}
\newcommand{\fp}{{\mathfrak p}}

\newcommand{\fs}{{\mathfrak s}}

\newcommand{\ch}{{\rm ch}}
\newcommand{\growth}{\mathop{\rm growth}}
\newcommand{\Lie}{\mathop{\rm Lie}}
\newcommand{\Par}{\mathop{\rm Par}}
\newcommand{\tr}{\mathop{\rm tr}}
\newcommand{\Res}{{\rm Res}}

\newtheorem*{corollary*}{Corollary}
\newtheorem*{remark*}{Remark}
\newtheorem*{remarks*}{Remarks}

\newcommand{\alphaparenlist}{
  \renewcommand{\theenumi}{\alph{enumi}}%
  \renewcommand{\labelenumi}{(\theenumi)}%
}

\makeatletter
\def\@maketitle{\newpage
 \begin{center}%
  \vskip 3em
 {\Large\bf \@title \par}%
  \vskip 1.5em
  {\normalsize
   \lineskip .5em
  \begin{tabular}[t]{c}\@author
   \end{tabular}\par}%
  \vskip 1em

 \end{center}%
 \par
 \vskip .5em}
\makeatother

\begin{document}

\title{Finite growth representations of infinite Lie conformal algebras}

\author{Carina Boyallian$^*$ \and Victor G. Kac$^\dag$ \and Jose 
I.~Liberati\thanks{
Famaf-CIEM, Ciudad Universitaria, (5000) Cordoba, Argentina,\newline
$<$boyallia@mate.uncor.edu$>$, $<$liberati@mate.uncor.edu$>$.
\newline
    $^\dag$ Department of Mathematics, M.I.T., 
   Cambridge, MA 02139, 
   $<$kac@math.mit.edu$>$.
}}

\maketitle

\begin{abstract} We classify all finite growth representations of all infinite 
rank subalgebras 
of the Lie conformal algebra $gc_1$ that contain a Virasoro subalgebra.
  
\end{abstract}

\section{Introduction} In this paper we study representation theory of some 
infinite rank subalgebras  of  the Lie conformal algebra $gc_1$ associated to 
the Lie algebra $\D$ of differential operators on the circle. Recall that 
$gc_1 = \CC [\p , x]$, with $\lambda$-bracket given by (see \cite{BKL}, 
\cite{K2})
\begin{displaymath}
	[a(\p , x)_\lambda b(\p , x)] = a(-\lambda ,\lambda +\p + x )b(\lambda 
+\p , x) - b(\lambda+\p, - \lambda +x) a(- \lambda, x)
\end{displaymath}

Also recall that the {\it Virasoro conformal algebra} (which is particularly 
important in physics) is defined as the free $\CC [\p ]$-module of rank 1 
generated by an element $L$, with $\lambda$-bracket defined by
\begin{displaymath}
[L_\lambda L]=(2\lambda +\p ) L,
\end{displaymath}
and extended to $\CC [\p ] L$ using sesquilinearity. Observe that all Virasoro 
subalgebras of $gc_1$ are generated by 
\begin{displaymath}
L=x +\alpha \p \, , \quad \alpha\in \CC.
\end{displaymath}

The complete list of infinite rank  proper subalgebras of $gc_1$ that contain a 
Virasoro subalgebra is (see \cite{BKL} and Remark 3.10 in \cite{DeK}):
\begin{eqnarray*}
	gc_{1,\, x} &=& x \, \CC[\p , x]  \\
	oc_1 &=& \{ a(\p , x) - a(\p , -\p -x) \, | \, a(\p ,x)\in \CC[\p , x] 
\}  \\
	spc_1 &=& \{ x [ a(\p , x) + a(\p , -\p -x)] \, | \, a(\p ,x)\in \CC[\p 
, x] \} 
\end{eqnarray*}
where the Virasoro element is $x+\alpha \p$, with $\alpha=0,\frac 1 2 , 0$ 
respectively. They are the most important $gc_1$-subalgebras from the point of 
view of physics.

In the present paper we classify all finite growth representations of all 
infinite rank conformal subalgebras of $gc_1$ that contain a Virasoro 
subalgebra.

This problem reduces to the  study of finite growth representations on the 
corresponding extended annihilation algebras, which are certain subalgebras of 
$\D$ (see \cite{K2}). The main tools used here are the recent results 
(\cite{AFMO},\cite{BL},\cite{KL},\cite{KR1} and \cite{KWY})  on the 
classification of quasifinite highest weight modules over the central extension 
of $\D$ and some of its important 
subalgebras.

The paper is organized as follows. In Section \ref{sec:2}, we describe the 
infinite rank Lie 
algebra $\widehat{g\ell}^{[m]}_{\infty}$ and its classical subalgebras, and 
discuss their
representation theory that will be needed. In sections 3,4,5, and 6, we obtain 
the classification 
of all finite growth representations of $gc_1$, $gc_{1,\,x}$, $oc_1$, and 
$spc_1$
respectively.

%
%

\section{Lie algebra $\widehat{g\ell}^{[m]}_{\infty}$ and its classical 
subalgebras}\label{sec:2}

\subsection{Lie algebra $\widehat{g\ell}^{[m]}_{\infty}$}

Denote by $\CC^{+\infty}$ the set of all sequences $\lambda =
(\lambda_1, \lambda_2, \ldots)$ for which all but a finite
number of $\lambda_i$'s are zero, and let $d (\lambda)$ denote the
number of non-zero $\lambda_i$'s and $|\lambda |$ denote their
sum.  Denote by $\Par^+$ the subset of $\CC^{+\infty}$ consisting 
of non-increasing sequences of (non-negative) integers.

Denote by $g\ell_{+\infty}$ the Lie algebra of all matrices
$(a_{ij})^{+\infty}_{i,j=1}$ with a finite number of non-zero
entries $a_{ij} \in \CC$.  Given $\lambda \in \CC^{+\infty}$,
there exists a unique irreducible $g\ell_{+\infty}$-module
$L^+(\lambda)$, also denoted by $L(g\ell_{+\infty};\lambda)$,  which admits a 
non-zero vector $v_{\lambda}$ such that
\begin{equation}
  \label{eq:1}
  E_{ij} v_{\lambda}=0 \/\/ \hbox{ for } \,\, i<j 
  \hbox{ and } E_{ii} v_{\lambda} =\lambda_i v_{\lambda} \, .
\end{equation}
Here and further $E_{ij}$ denotes, as usual, the matrix whose
$(i,j)$-entry is $1$ and all other entries are $0$.  Each
$L^+(\lambda)$ has a unique $\ZZ_+$-gradation.  $\displaystyle{L^+(\lambda)
=\oplus_{j \in \ZZ_+}}$ $L^+(\lambda)_j$, called its \emph{principal gradation}, 
which satisfies the properties 
\begin{displaymath}
  L^+(\lambda)_0 = \CC v_{\lambda} \, , \quad
  E_{ij} L^+(\lambda)_k \subset 
  L^+ (\lambda)_{k+i-j} \, .
\end{displaymath}
Since $\lambda \in \CC^{+\infty}$, it is easy to see that $\dim
L^+(\lambda)_j< \infty$, hence we can define the $q$-character
\begin{displaymath}
  \ch_q L^+ (\lambda) =
  \sum_{j \in \ZZ_+} (\dim L^+ (\lambda)_j)q^j  \, .
\end{displaymath}

For $\lambda \in \Par^+$, let $d=d (\lambda)$ and $\bar{\lambda}
= (\lambda_1 , \ldots , \lambda_d)$.  Let $g\ell_d$ be the Lie
algebra of all $d \times d$ matrices $(a_{ij})^d_{i,j=1}$; it may 
be viewed as a subalgebra of $g\ell_{+\infty}$ in a natural way.
Denote by $\bar{L}^+(\bar{\lambda})$ the (irreducible)
$g\ell_d$-submodule of $L^+(\lambda)$ generated by
$v_{\lambda}$.  It is, of course, isomorphic to the
finite-dimensional irreducible $g\ell_d$-module associated to
$\bar{\lambda}$, so that its $q$-character is a (well-known)
polynomial in $q$.

\begin{lemma}
  \label{lem:1}
Let $\lambda \in \Par^+, d=d(\lambda)$.  Then 
\begin{displaymath}
  \ch_q L^+ (\lambda) = ch_q \bar{L}^+ (\bar{\lambda})/
  \prod^d_{j=1} (1-q^j)^{\lambda_{d-j+1}}_{q} \, ,
\end{displaymath}
where  $(1-a)^m_q = (1-a) (1-qa) \ldots (1-q^{m-1}a)$.

\end{lemma}

\begin{proof}
  Recall the well-known formula (see \cite{K1}):
  \begin{equation}
    \label{eq:2}
    \ch_q L^+ (\lambda) = \prod_{\alpha >0}
    (1-q^{\langle \lambda + \rho , \alpha \rangle})/
    (1-q^{\langle \rho , \alpha \rangle}) \, .    
  \end{equation}
Here the product is taken over the set of all positive coroots of 
$g\ell_{+\infty}$, which are all elements $E_{ii}-E_{jj}$ with
$i<j$, $\langle \lambda ,E_{ii} \rangle =\lambda_i$ and $\langle
\rho , E_{ii} \rangle =-i$.  Of course, a similar formula holds
for $\ch_q \bar{L}^+ (\bar{\lambda})$; it is a part of the
product~(\ref{eq:2}) corresponding to $i<j \leq d$.

It is also clear that the factors of (\ref{eq:2}) corresponding
to $d<i<j$ are equal to $1$, and it is easy to see that the
product over all pairs $i,j$ with fixed $i \leq d$ and all $j>d$
is equal to $1/(1-q^{d-i+1})^{\lambda_i}_q$.

\end{proof}

Recall that, given a vector space $V$ with an increasing
filtration by finite-dimensional subspaces $V_{[j]}$, the
\emph{growth} of $V$ is defined by
\begin{displaymath}
  \growth V = \overline{\lim}_{j \to +\infty}
  (\log \dim V_{[j]})/\log j \, .
\end{displaymath}
We define the growth of $L^+(\lambda)$ using its filtration
$L^+(\lambda)_{[j]}=\oplus_{i \leq j} L^+(\lambda)_i$ associated to 
the principal gradation.

\begin{theorem}
  \label{th:1}
  \alphaparenlist
  \begin{enumerate}
  \item 
\label{th:1a}
If $\lambda \in \Par^+$, then
\begin{displaymath}
   \growth  L^+ (\lambda) = |\lambda | \, .
\end{displaymath}

\item 
\label{th:1b}
If $\lambda \in \CC^{+\infty} \backslash \Par^+$, then $\growth
L^+ (\lambda ) =\infty$.

  \end{enumerate}
\end{theorem}

\begin{proof}
  It follows from Lemma~\ref{lem:1} that for $\lambda \in \Par^+$, 
  $\growth L^+(\lambda)$ is equal to the growth of the
  polynomial algebra on generators of degree $1,2, \ldots
  ,\lambda_s$; $2,3, \ldots ,\lambda_{s-1}+1$; $\ldots$; $s,s+1
  ,\ldots , \lambda_1 +s-1$.  The total number of these
  generators is $|\lambda |$, and since growth of a polynomial
  algebra is independent of the degrees of generators, (a)~is
  proved.

Let now $\lambda \in \CC^{+\infty} \backslash \Par^+$.  Then
$\lambda_k -\lambda_{k+1} \notin \ZZ_+$ for some $k$.  But then $
E^N_{k+1,k} v_{\lambda }\neq 0$ for each $N \in \ZZ_+$.  Looking
at the subalgebra of $g\ell_{+\infty}$ spanned by all $E_{ij}$
with $i,j \geq k+1$, we conclude from~(a) that
\begin{displaymath}
  \growth L^+ (\lambda ) \geq N+\sum_{i \geq k+1} \lambda_i \, .
\end{displaymath}
This proves (b).
\end{proof}

In a similar fashion one may consider the Lie algebra
$g\ell_{-\infty}$ of all matrices $(a_{ij})^{-\infty}_{i,j=0}$
with a finite number of non-zero entries and the irreducible
$g\ell_{-\infty}$-modules $L^-(\lambda)$, also denoted by $L(g\ell_{-\infty}; 
\lambda)$, parameterized by the set 
$\CC^{-\infty}$ of sequences $\mu = (\ldots , \mu_{-1}, \mu_0)$
with finitely many non-zero members.  Results similar to
Lemma~\ref{lem:1} and Theorem~\ref{th:1} hold for the subset
$\Par^- \subset \CC^{-\infty}$ consisting of non-increasing
sequences of (non-positive) integers.

Let $\tilde{g\ell}_{\infty}$ denote the Lie algebra of all
matrices $(a_{ij})_{i,j\in \ZZ}$ such that $a_{ij}=0$ if $|i-j|
\gg 0$.  Denote by $\tilde{g\ell}_{+\infty}$
(resp. $\tilde{g\ell}_{-\infty}$) the subalgebra of
$\tilde{g\ell}_{\infty}$ consisting of matrices with $a_{ij}=0$
for $i$ or $j \leq 0$ (resp. $i$ or $j>0$).  Note that these two
subalgebras commute and that $\tilde{g\ell}_{\pm \infty}$
contains $g\ell_{\pm \infty}$ as a subalgebra.  Note also that
the $g\ell_{\pm \infty}$-modules $L^{\pm}(\lambda)$ extend
uniquely to $\tilde{g\ell}_{\pm\infty}$.

The Lie algebra $\tilde{g\ell}_{\infty}$ has a well-known central 
extension $\widehat{g\ell}_{\infty} =\tilde{g\ell} +\CC C$ by $\CC$
defined by the cocycle
\begin{equation}\label{eq:co}
  \alpha (A,B) = \tr [J,A]B \, , \,\, \hbox{ where }
  J=\sum_{i \leq 0} E_{ii} \, .
\end{equation}
The restriction of this cocycle to $\tilde{g\ell}_{+\infty}$ and
to $\tilde{g\ell}_{-\infty}$ is zero.

We will also need briefly the Lie algebra
$\widehat{g\ell}^{[m]}_{\infty}$ defined for each $m \in \ZZ_+$ by
replacing $\CC$ by $R_m=\CC[u]/(u^{m+1})$. That is, 
$\widehat{g\ell}^{[m]}_{\infty}=\tilde{g\ell}^{[m]}_{\infty}\oplus R_m$ is the 
central extension of $\tilde{g\ell}^{[m]}_{\infty}$ by the 2-cocycle 
(\ref{eq:co}) with values in $R_m$, where $\tilde{g\ell}^{[m]}_{\infty}$ is the 
Lie algebra of infinite matrices with finitely many nonzero diagonals with 
entries in $R_m$.

The principal $\ZZ$-gradation of all the above Lie algebras is
defined by letting
\begin{equation}\label{eq:gradation}
  \deg E_{ij}=i-j
\end{equation}
(in the case of $\widehat{g\ell}^{[m]}_{\infty}$ we also let $\deg R_m=0$). This 
gives us a triangular decomposition
\begin{displaymath}
   \widehat{g\ell}^{[m]}_{\infty}=(\widehat{g\ell}^{[m]}_{\infty})_+ \oplus 
(\widehat{g\ell}^{[m]}_{\infty})_0 \oplus (\widehat{g\ell}^{[m]}_{\infty})_-,
\end{displaymath}
where
\begin{displaymath}
  (\widehat{g\ell}^{[m]}_{\infty})_\pm = \oplus_{j\in \NN} 
(\widehat{g\ell}^{[m]}_{\infty})_{\pm j }
\end{displaymath}

The Lie algebra $\widehat{g\ell}_{\infty}$ has a family of modules
$L(\widehat{g\ell}_{\infty};\lambda,c)$, parameterized by $\lambda \in 
\CC^{\infty} = \{
(\lambda_i)_{i \in \ZZ}|$ all but finitely many of $\lambda_i$
are $0\}$ and $c \in \CC$, defined by (\ref{eq:1}) and $C
v_{\lambda} = cv_{\lambda}$.  Similarly,
$\widehat{g\ell}^{[m]}_{\infty}$ has a family of modules
 $L(\widehat{g\ell}^{[m]}_{\infty}; \vec{\lambda},\vec{c})$, where 
$\vec{\lambda} \in
(\CC^{\infty})^{m+1}$, $\vec{c} \in \CC^{m+1}$, defined in a
similar fashion. That is,  the highest weight 
$\widehat{g\ell}^{[m]}_{\infty}$-module $L(\widehat{g\ell}^{[m]}_{\infty}; 
\Lambda)$, with highest weight $\Lambda\in (\widehat{g\ell}^{[m]}_{\infty})_0^*$ 
 that is determined by its {\it labels} $\vec\lambda_i^{(j)}=\Lambda(u^j 
E_{ii})$ and the {\it central charges} $\vec c_j=\Lambda(u^j)$.

The gradation (\ref{eq:gradation}) is obviously consistent with the principal
gradation of $L^{\pm}(\lambda)$ and of $L(\widehat{g\ell}_{\infty};\lambda ,c)$.

%
%

\subsection{Lie algebras $b_{\infty} ^{[m]}$ and $d_{\infty} ^{[m]}$ }

The Lie algebra $\tilde{g\ell}_{\infty}^{[m]}$ acts on the vector space $R_m[t, 
t^{-1}]$ via the usual formula
\begin{displaymath}
     E_{ij}v_k =\delta_{j,k}v_i \, ,
\end{displaymath}
where $v_i=t^{-i},\, i\in\ZZ$ is an $R_m$-basis. 
Now consider the following $\CC$-bilinear forms on this space:
\begin{eqnarray*}
   B(u^m v_i,u^nv_j)&=& u^{m}(-u^{n})\delta_{i,-j}  \\
   D(u^m v_i, u^n v_j)&=& u^{m}(-u^{n})\delta_{i,1-j}  \, .
\end{eqnarray*}
Denote by $\bar b_{\infty}^{[m]}$ (resp. $\bar d_{\infty}^{[m]}$) the Lie 
subalgebra of $\tilde{g\ell}_{\infty}^{[m]}$ which preserves the bilinear form 
$B$ (resp. $D$). We have
\begin{eqnarray*}
   \bar b_{\infty}^{[m]}&=& \{ (a_{ij}(u))_{i,j\in\ZZ} \in 
\tilde{g\ell}_{\infty} ^{[m]}\, | \, a_{ij}(u)=-a_{-j,-i}(-u)\, \} \\
   \bar d_{\infty}^{[m]}&=&  \{ (a_{ij}(u))_{i,j\in\ZZ} \in 
\tilde{g\ell}_{\infty}^{[m]} \, | \, a_{ij}(u)=-a_{1-j,1-i}(u)\, \}  \, .
\end{eqnarray*}
Denote by $b_{\infty} ^{[m]}=\bar b_{\infty} ^{[m]}\oplus R_m$ (resp. 
$d_{\infty} ^{[m]}=\bar d_{\infty} ^{[m]}\oplus R_m$) the central extension of 
$\bar b_{\infty} ^{[m]}$ (resp. $\bar d_{\infty} ^{[m]}$) given by the 2-cocycle 
 defined in $\tilde g\ell_\infty ^{[m]}$. Both subalgebras inherits form 
$\widehat {g\ell}_\infty ^{[m]}$ the principal $\ZZ$-gradation and the 
triangular decomposition, (see \cite{KWY} and \cite{K1} for notation)
\begin{eqnarray*}
    b_{\infty}^{[m]}&=&\oplus_{j\in\ZZ} (b_{\infty}^{[m]})_j\qquad 
b_{\infty}^{[m]}=(b_{\infty}^{[m]})_+\oplus (b_{\infty}^{[m]})_0\oplus 
(b_{\infty}^{[m]})_-\\
    d_{\infty}^{[m]}&=&\oplus_{j\in\ZZ} (d_{\infty}^{[m]})_j\qquad 
d_{\infty}^{[m]}=(d_{\infty}^{[m]})_+\oplus (d_{\infty}^{[m]})_0\oplus 
(d_{\infty}^{[m]})_- \, .
\end{eqnarray*}

In particular when $m=0$, we have the usual Lie subalgebras of 
$\widehat{g\ell}_\infty$, denoted by $b_\infty$ (resp. $d_\infty$).

Denote by $L( b_{\infty}^{[m]}; \lambda)$ (resp. $L( d_{\infty}^{[m]}; 
\lambda)$) the highest weight module over $b_{\infty}^{[m]}$ (resp. 
$d_{\infty}^{[m]}$) with highest weight $\lambda\in (b_{\infty}^{[m]})^*_0$ 
(resp. $\lambda\in (d_{\infty}^{[m]})^*_0$)   parametrized by $^b\vec{\lambda} 
\in
(\CC^{\infty})^{m+1}$, $\vec{c} \in \CC^{m+1}$, with
\begin{eqnarray*}
\vec c_i &= & \lambda(u^i), \\
  \quad\quad \ ^b\vec\lambda_j^{(i)}&=&\lambda(u^i\,E_{j,j}-(-u)^i\, E_{-j,-j}),
\end{eqnarray*}
(resp. $^d\vec{\lambda} \in
(\CC^{\infty})^{m+1}$, with $^d\vec\lambda_j^{(i)}=\lambda(u^i\,E_{j,j}-(-u)^i\, 
E_{1-j,1-j})$). The superscripts $b$ and $d$ here mean B and D type 
respectively. The $^b\vec\lambda_j^{(i)}$ (resp. $^d\vec\lambda_j^{(i)}$) are 
called the labels and $\vec c_j$ the central charges of $L( 
b_\infty^{[m]};\lambda)$ (resp. $L( d_\infty^{[m]};\lambda)$).

All these  modules will appear in Section~\ref{sec:oc1}. Now, we are interested 
in  representation theory of $b_{\infty}$. 

 The set of simple coroots of $b_\infty$, can be described as follows (cf. 
\cite{KWY}): 
\begin{eqnarray*}
    \Pi\,\check{ } =\{{\alpha_0}\check{ } &=&  2(E_{-1,-1}-E_{1,1})+2C,  \\
  \alpha_i\check{ }  &=& E_{i,i}-E_{i+1,i+1}-E_{-i,-i}+E_{-1-i,-1-i} 
\,,\,i\in\NN\}  \, .
\end{eqnarray*}
The set of roots:
\begin{displaymath}
   \Delta=\{\pm\varepsilon_0,\, \pm\varepsilon_i,\, \pm\varepsilon_i 
\pm\varepsilon_j,\, i\neq j,\,i,j\in\NN\}  \, .
\end{displaymath}
The set of positive coroots is:
\begin{eqnarray*}
    \Delta_+\check{ }\,&=&\{\alpha_i\check{ }+ \alpha_{i+1}\check{ }+\cdots 
+\alpha_j\check{ }\, ,\, 0\leq i\leq j\}\, \\ 
&\cup& \{\alpha_0\check{ }+ 2\alpha_{ 1}\check{ }+\cdots +2\alpha_{i}\check{ 
}+\alpha_{i+1}\check{ }  +\cdots +\alpha_{j-1}\check{ }\, ,\, 1\leq i< j\}  \, .
\end{eqnarray*}
The set of simple roots:
\begin{displaymath}
   \Pi=\{\alpha_0=-\varepsilon_1,\, \alpha_i=\varepsilon_i-\varepsilon_{i+1}  ,  
 \,i \in\NN\}  \, .
\end{displaymath}
Here $\varepsilon_i $ are viewed restricted to the restricted dual of the Cartan 
subalgebra of $b_\infty$, so that $\varepsilon_i=-\varepsilon_{-i}$. Given 
$\lambda\in(b_\infty)_0^*$, the labels and central charge are simply (in this 
case, we skip the superscript b)
\begin{displaymath}
   \lambda_i=\lambda(E_{i,i}-E_{-i,-i}), \, i>0\,, \quad c=\lambda(C) .
\end{displaymath}
So that  $\lambda(\alpha_0\check{ }\,)=2c-2\lambda_1$ and 
$\lambda(\alpha_i\check{ }\,)= \lambda_i-\lambda_{i+1}$ for $i\in\NN$. Denote by 
$\Lambda_i$ the i-th fundamental weight of $b_\infty$, namely 
$\Lambda_i(\alpha_j\check{ }\,)=\delta_{i,j}$.

Let $P_+=\{\lambda\in (b_\infty)_0^*\,|\, \langle \lambda,\alpha_i {\check{ }} 
\, \rangle\in\ZZ_+,\hbox{ for all }i\in\ZZ_+\}$ denote the set of dominant 
integral weights of $b_\infty$. Given $\lambda \in P_+$, we have 
$\lambda=\Lambda_{n_1}+\Lambda_{n_2}+\cdots+\Lambda_{n_k}+h\Lambda_0$, $n_1\geq 
n_2\geq\cdots\geq n_k\geq 1$, $h\in\ZZ_+$, and the module $L( b_\infty 
;\lambda)$ has central charge $c=k+h/2$. Observe that the conjugate of the Young 
diagram corresponding to the partition $(n_1,n_2,\cdots,n_k)$ is 
$(\lambda_1,\cdots,\lambda_{n_1})$, and $\lambda_i=0$ for $i>n_1$. Note that 
$n_1=n_1 (\lambda)=\max \{i\in\NN \,|\, \langle \lambda,\alpha_i {\check{ 
}}\rangle\neq 0\}$. Observe that $\fs\fo(2n_1+1)$ may be viewed as a subalgebra 
of $b_\infty$  in a natural way, whose set of simple roots is $\{-\varepsilon_1, 
\varepsilon_1 - \varepsilon_2, \dots, \varepsilon_{n_1-1} - 
\varepsilon_{n_1}\}$. Denote by    $\bar{\lambda}$ the dominant integral weight 
of $\fs\fo(2n_1+1)$ g!
iven by $\bar{\lambda}( 2(E_{-1,-1}-E_{1,1}))= 2(c-\lambda_1 )$ and  
$\bar{\lambda}( E_{i,i}-E_{i+1,i+1}-E_{-i,-i}+E_{-1-i,-1-i})=\lambda_{i} - 
\lambda_{i+1}$ for $1\leq i < n_1 $. Denote by $\bar{L} (\bar{\lambda})$ the 
(irreducible)
 $\fs\fo(2n_1+1)$-submodule of $L (b_\infty ;\lambda)$ generated by
its highest weight vector.  It is, of course, isomorphic to the
finite-dimensional irreducible  $\fs\fo(2n_1+1)$-module associated to
$\bar{\lambda}$, so that its $q$-character is a (well-known)
polynomial in $q$.

\begin{lemma}
  \label{lem:5}
Let $\lambda \in P_+, n_1=n_1(\lambda)$.  Then 
\begin{eqnarray*}
  \ch_q L  (b_\infty;\lambda) = ch_q \bar{L}  (\bar{\lambda}) . 
  \prod^{n_1}_{j=1} \frac{1}{(1-q^j)^{\lambda_{n_1-j+1}}_{q} } &.&
  \prod^{n_1}_{i=1} \frac{1}{(1-q^{n_1+i})^{2c - \lambda_i}_{q} }\, \\
 &\times& \prod_{n_1\leq i} \frac{1}{(1-q^{2i+1})^{2c}_{q} } ,
\end{eqnarray*}
where $ (1-a)^m_q = (1-a) (1-qa) \ldots (1-q^{m-1}a)$.

\end{lemma}

\begin{proof} The proof is completely similar to the one of Lemma~{\ref{lem:1}}, 
using the data introduced above (cf. proof of Proposition 1.1 in \cite{KWY}).

\end{proof}

\begin{theorem}
  \label{th:5}  All non-trivial modules $L( b_{\infty}^{[m]}; \lambda)$ have 
infinite growth.

\end{theorem}

\begin{proof} It is enough to consider the case $m=0$. Given $\lambda\in 
(b_\infty)_0^*$, we look 
at the subalgebra of $b_\infty$ isomorphic to $g\ell_{+\infty}$ spanned by all 
$E_{i,j}- E_{-j,-i}$
with $i,j \geq 1$, and by Theorem~{\ref{th:1}}, we conclude that $L(b_\infty ; 
\lambda)$ has infinite growth if $\lambda_i -\lambda_{i+1}\notin\ZZ_+$ for some 
$i\geq 1$.

Let us assume $\lambda_i -\lambda_{i+1}\in\ZZ_+$ for all $i\geq 1$. If 
$2c-2\lambda_1 \notin\ZZ_+$, then  $
(E_{ 1,0}-E_{0,-1})^N v_{\lambda }\neq 0$ for each $N \in \ZZ_+$. 
Looking at the subalgebra of $b_\infty$ isomorphic to $g\ell_{+\infty}$ 
previously defined, we conclude from Theorem~{\ref{th:1}} that 
\begin{displaymath}
  \growth L  (b_\infty ;\lambda ) \geq N+\sum_{i \geq  1} \lambda_i \, .
\end{displaymath}
If $2c-2\lambda_1 \in\ZZ_+$, then by the same argument in the proof of 
Theorem~{\ref{th:1}}, and looking at the last factor in Lemma~\ref{lem:5}, we 
conclude that $L(b_\infty ; \lambda)$ has the same growth as the polynomial 
algebra in infinitely many generators, finishing the proof.

\end{proof}

%
%

\subsection{Lie algebra $c_{\infty} ^{[m]}$}

As before, we consider the vector space $R_m[t, t^{-1}]$, and take the 
$R_m$-basis $v_i=t^{-i},\, i\in\ZZ$.
Denote by $\bar c_{\infty}^{[m]}$  the Lie subalgebra of 
$\tilde{g\ell}_{\infty}^{[m]}$ which preserves the bilinear form:
\begin{equation}
   C(u^m v_i,u^nv_j)= u^{m}(-u^{n}) (-1)^i \, \delta_{i,1-j}   \, .
\end{equation}
We have
\begin{displaymath}
   \bar c_{\infty}^{[m]}= \{ (a_{ij}(u))_{i,j\in\ZZ} \in \tilde{g\ell}_{\infty} 
^{[m]}\, | \, a_{ij}(u)=(-1)^{i+j+1}a_{1-j,1-i}(-u)\, \} \, .\end{displaymath}
Denote by $c_{\infty} ^{[m]}=\bar c_{\infty} ^{[m]}\oplus R_m$   the central 
extension of $\bar c_{\infty} ^{[m]}$   given by the 2-cocycle  defined in 
$\tilde g\ell_\infty ^{[m]}$. This subalgebra inherits form $\widehat 
{g\ell}_\infty ^{[m]}$ the principal $\ZZ$-gradation and the triangular 
decomposition, (see \cite{KWY} and \cite{K1} for notation)
\begin{displaymath}
    c_{\infty}^{[m]}=\oplus_{j\in\ZZ} (c_{\infty}^{[m]})_j\qquad 
c_{\infty}^{[m]}=(c_{\infty}^{[m]})_+\oplus (c_{\infty}^{[m]})_0\oplus 
(c_{\infty}^{[m]})_- \, .
\end{displaymath}

In particular when $m=0$, we have the usual Lie subalgebra of 
$\widehat{g\ell}_\infty$, denoted by $c_\infty$ (see \cite{K1}).

Denote by $L( c_{\infty}^{[m]}; \lambda)$ the highest weight module over 
$c_{\infty}^{[m]}$  with highest weight $\lambda\in (c_{\infty}^{[m]})^*_0$  
parametrized by its labels $^c\vec{\lambda} \in
(\CC^{\infty})^{m+1}$ and central charges  $\vec{c} \in \CC^{m+1}$, with
\begin{eqnarray*}
\vec c_i &= & \lambda(u^i), \\
  \quad\quad \ ^c\vec\lambda_j^{(i)}&=&\lambda(u^i\,E_{j,j}-(-u)^i\, 
E_{1-j,1-j}).
\end{eqnarray*}

Now, we are interested in  representation theory of $c_{\infty}$. 

 The set of simple coroots of $c_\infty$, denoted by $\Pi\,\check{ }$, can be 
described as follows (cf. \cite{KWY}): 
\begin{eqnarray*}
    \Pi\,\check{ }=\{\alpha_0\check{ }&=&  E_{0,0}-E_{1,1}+C,  \\
  \alpha_i\check{ }  &=& E_{i,i}-E_{i+1,i+1}+E_{-i,-i}-E_{1-i,1-i} 
\,,\,i\in\NN\}  \, .
\end{eqnarray*}
The set of roots:
\begin{displaymath}
   \Delta=\{ \pm 2 \varepsilon_i,\, \pm\varepsilon_i \pm\varepsilon_j,\, i\neq 
j,\,i,j\in\NN\}  \, .
\end{displaymath}
The set of positive coroots is:
\begin{eqnarray*}
    \Delta_+\,\check{ }&=&\{\alpha_i\check{ }+ \alpha_{i+1}\check{ }+\cdots 
+\alpha_j\check{ }\, ,\, 0\leq i\leq j\}\, \\ 
&\cup& \{2\alpha_0\check{ }+ 2\alpha_{ 1}\check{ }+\cdots +2\alpha_{i}\check{ 
}+\alpha_{i+1}\check{ }   +\cdots +\alpha_{j}\check{ }\, ,\, 0\leq i< j\}  \, .
\end{eqnarray*}
The set of simple roots:
\begin{displaymath}
   \Pi=\{\alpha_0=-2\varepsilon_1,\, \alpha_i=\varepsilon_i-\varepsilon_{i+1}  , 
  \,i \in\NN\}  \, .
\end{displaymath}
Here $\varepsilon_i $ are viewed restricted to the restricted dual of the Cartan 
subalgebra of $c_\infty$, so that $\varepsilon_i=-\varepsilon_{1-i}$. Given 
$\lambda\in(c_\infty)_0^*$, the labels and central charge are simply (in this 
case, we skip the superscript c):
\begin{displaymath}
   \lambda_j =\lambda(E_{j,j}- E_{1-j,1-j}),\quad j\in\NN,     	\quad 
c=\lambda(C).
\end{displaymath} 
So that  $\lambda(\alpha_0\check{ }\,)=-\lambda_1+c$ and 
$\lambda(\alpha_i\check{ }\,)= \lambda_i-\lambda_{i+1}$ for $i\in\NN$. Denote by 
$\Lambda_i$ the i-th fundamental weight of $b_\infty$, namely 
$\Lambda_i(\alpha_j\check{ }\,)=\delta_{i,j}$.

Let $P_+=\{\lambda\in (c_\infty)_0^*\,|\, \langle \lambda,\alpha_i {\check{ }} 
\, \rangle\in\ZZ_+,\hbox{ for all }i\in\ZZ_+\}$ denote the set of dominant 
integral weights of $c_\infty$. Given $\lambda \in P_+$, we have 
$\lambda=\Lambda_{n_1}+\Lambda_{n_2}+\cdots+\Lambda_{n_k}+h\Lambda_0$, $n_1\geq 
n_2\geq\cdots\geq n_k\geq 1$, $h\in\ZZ_+$, and the module $L( c_\infty,\lambda)$ 
has central charge $c=k+h$. Observe that the conjugate of the Young diagram 
corresponding to the partition $(n_1,n_2,\cdots,n_k)$ is 
$(\lambda_1,\cdots,\lambda_{n_1})$, and $\lambda_i=0$ for $i>n_1$. Note that 
$n_1=n_1 (\lambda)=\max \{i\in\NN \,|\, \langle \lambda,\alpha_i {\check{ 
}}\rangle\neq 0\}$. Observe that $\fs\fp(2n_1)$ may be viewed as a subalgebra of 
$c_\infty$  in a natural way, whose set of simple roots is $\{-2\varepsilon_1, 
\varepsilon_1 - \varepsilon_2, \dots, \varepsilon_{n_1-1} - 
\varepsilon_{n_1}\}$. Denote by    $\bar{\lambda}$ the dominant integral weight 
of $\fs\fp(2n_1)$ given b!
y $\bar{\lambda}( (E_{0,0}-E_{1,1}))= c-\lambda_1 $ and  $\bar{\lambda}( 
E_{i,i}-E_{i+1,i+1}+E_{-i,-i}-E_{1-i,1-i})=\lambda_{i} - \lambda_{i+1}$ for 
$1\leq i < n_1 $. Denote by $\bar{L} (\bar{\lambda})$ the (irreducible)
 $\fs\fp(2n_1)$-submodule of $L (c_\infty,\lambda)$ generated by
its highest weight vector.  It is, of course, isomorphic to the
finite-dimensional irreducible  $\fs\fp(2n_1)$-module associated to
$\bar{\lambda}$, so that its $q$-character is a (well-known)
polynomial in $q$.

\begin{lemma}
  \label{lem:4.1}
Let $\lambda \in P_+, n_1=n_1(\lambda)$.  Then 
\begin{eqnarray*}
  \ch_q L  (c_\infty;\lambda) = ch_q \bar{L}  (\bar{\lambda}) &.& 
  \prod^{n_1}_{i=1} \frac{1}{(1-q^j)^{\lambda_{n_1-j+1}}_{q} } .
  \prod^{n_1}_{i=1} \frac{1}{(1-q^{n_1+i+3})^{2c - \lambda_{i+1}}_{q} }\, \\
 &\times& \frac{1}{(1-q^{n_1+1})^c_q} . \prod_{n_1\leq i} 
\frac{1}{(1-q^{2i+3})^{2c}_{q} } ,
\end{eqnarray*}
where $ (1-a)^m_q = (1-a) (1-qa) \ldots (1-q^{m-1}a)$.

\end{lemma}

\begin{proof} The proof is completely similar to the one of Lemma~{\ref{lem:1}}, 
using the data introduced above (cf. proof of Proposition 1.1 in \cite{KWY}).

\end{proof}

\begin{theorem}
  \label{th:4.2}  All non-trivial modules $L( c_{\infty}^{[m]}; \lambda)$ have 
infinite growth.

\end{theorem}

\begin{proof} It is enough to consider the case $m=0$. Given $\lambda\in 
(c_\infty)_0^*$, we look 
at the subalgebra of $c_\infty$ isomorphic to $g\ell_{+\infty}$ spanned by all 
$E_{i,j}- E_{1-j,1-i}$
with $i,j \geq 1$, and by Theorem~{\ref{th:1}}, we conclude that $L(c_\infty , 
\lambda)$ has infinite growth if $\lambda_i -\lambda_{i+1}\notin\ZZ_+$ for some 
$i\geq 1$.

Let us assume $\lambda_i -\lambda_{i+1}\in\ZZ_+$ for all $i\geq 1$. If 
$c-\lambda_1\notin\ZZ_+$, then  $
(E_{ 1,0})^N v_{\lambda }\neq 0$ for each $N \in \ZZ_+$. 
Looking at the subalgebra of $c_\infty$ isomorphic to $g\ell_{+\infty}$ 
previously defined, we conclude from Theorem~{\ref{th:1}} that 
\begin{displaymath}
  \growth L  (c_\infty,\lambda ) \geq N+\sum_{i \geq  1} \lambda_i \, .
\end{displaymath}
If $c-\lambda_1 \in\ZZ_+$, then by the same argument in the proof of 
Theorem~{\ref{th:1}}, and looking at the last factor in Lemma~\ref{lem:5}, we 
conclude that $L(c_\infty , \lambda)$ has the same growth as the polynomial 
algebra in infinitely many generators, finishing the proof.

\end{proof}

%
%

\section{Irreducible finite growth $gc_1$-modules}\label{sec:3}

We are interested in representation theory of the Lie algebra
$\D^-$ of regular differential operators on $\CC$.  It consists of 
linear combinations of differential operators of the form $f(t)
\left( \frac{d}{dt} \right)^m$, where $f$ is a polynomial and $m
\in \ZZ_+$.  In particular, $D=t \frac{d}{dt} \in \D^-$.  The
principal $\ZZ$-gradation $\D^- = \oplus_{j \in \ZZ} \D^-_j$ is
defined by letting
\begin{displaymath}
  \deg t = -1 \, \,\, \deg \frac{d}{dt} =1 \, .
\end{displaymath}

Given a sequence of complex numbers $\Delta = (\Delta_0 ,
\Delta_1 , \ldots)$ we define the highest weight module $L(\Delta 
; \D^-)$ over $\D^-$ as the (unique) irreducible module that has
a non-zero vector $v_{\Delta}$ with the following properties:
\begin{displaymath}
  \D^-_j v_{\Delta} =0 \hbox{ for } j<0 , \,\,
    D^n v_{\Delta} = \Delta_n v_{\Delta}
    \hbox{ for }n \in \ZZ_+ \, .
\end{displaymath}
The principal gradation of $\D^-$ induces the principal gradation $L (\Delta ; 
\D^-)
=\oplus_{j \in \ZZ_+} L_j$ such that $L_0 =\CC v_{\Delta}$.  The
module $L (\Delta ; \D^-)$ is called \emph{quasifinite} if $\dim
L_j <\infty$ for all $j \in \ZZ_+$.

Quasifinite modules over $\D^-$ can be constructed as follows.
Consider the natural action of $\D^-$ on $\CC [t,t^{-1}]$ and
choose the basis $v_j=t^{-j} (j \in \ZZ)$ of $\CC [t,t^{-1}]$.
This gives an embedding of $\D^-$ in $\tilde{g\ell}_{\infty}$.
Since $\CC [t]$ is $\D^-$-invariant, we get $\D^-$-modules $\CC
[t,t^{-1}]/\CC [t]$ and $\CC [t]$, which gives us an embedding of 
$\D^-$ in $\tilde{g\ell}_{+\infty}$ and
$\tilde{g\ell}_{-\infty}$, respectively, hence an embedding of
$\D^-$ in $\tilde{g}\ell_{+\infty} \oplus
\tilde{g}\ell_{-\infty}$.  All these embeddings respect the
principal gradations.  Now take $\lambda^{\pm} \in \CC^{\pm
  \infty}$ and consider the $\tilde{g\ell}_{+\infty} \oplus
\tilde{g\ell}_{-\infty}$-module $L^+ (\lambda^+) \otimes L^-
(\lambda^-)$.  The same argument as in \cite{KR1}, gives us the
following.

\begin{lemma}
  \label{lem:2}
When restricted to $\D^-$, the module $L^+ (\lambda^+) \otimes
L^- (\lambda^-)$ remains irreducible.

\end{lemma}

It follows immediately that $L^+ (\lambda^+) \otimes L^-
(\lambda^-)$ is an irreducible highest weight module over $\D^-$, 
which is obviously quasifinite.  It is easy to see that we have:
\begin{eqnarray*}
  \Delta_n &=& \sum_{j \geq 1} (-j)^n \lambda^+_j
  + \sum_{j \leq 0} (-j)^n \lambda^-_j \, ,\\
\noalign{\hbox{so that}}\\
\Delta (x): &=& \sum_{n \geq 0} \Delta_n x^n /n! 
   = \sum_{j \geq 1}\lambda^+_j e^{-jx} 
   + \sum_{j \leq 0} \lambda^-_j e^{-jx} \, .
\end{eqnarray*}
It is also clear that for $\lambda^{\pm} \in \Par^{\pm}$ we have
(cf.~Theorem~1a):
\begin{displaymath}
  \growth L^+ (\lambda^+) \otimes L^- (\lambda^-)
  = |\lambda^+| + |\lambda^-| \, .
\end{displaymath}

We shall prove the following theorem.

\begin{theorem}
  \label{th:2}
  The $\D^-$-modules $L^+ (\lambda^+) \otimes L^- (\lambda^-)$,
  where $\lambda^{\pm} \in \Par^{\pm}$, 
  exhaust all quasifinite irreducible highest weight
  $\D^-$-modules that have finite growth.
\end{theorem}

Let $\D$ denote the Lie algebra of all regular differential
operators on $\CC^*$.  The Lie algebra $\D$ is the linear span of 
differential operators $\displaystyle{f(t) \left(
    \frac{d}{dt}\right)^k}$, where $f(t) \in \CC [t,t^{-1}]$ and
$k \in \ZZ_+$, or, equivalently of operators $t^k f(D)$, where
$f(D) \in \CC [D]$ and $k \in \ZZ$.  Obviously, $\D^-$ is a subalgebra of $\D$,
and the principal gradation extends from $\D^-$ to $\D$ in the
obvious way.

The basic idea of the proof of Theorem~\ref{th:2} is the same as
in \cite{}:  to reduce the problem to the well developed (in
\cite{KR1}) representation theory of the universal central extension 
$\hat{\D}$ of $\D$.  Recall that the central extension $\hat{\D}
= \D +\CC C$ is defined by the cocycle \cite{KR1}
\begin{equation}\label{eq:cocycle-d}
  \Psi \left( f (t)\left( \frac{d}{dt}\right)^m , \, 
    g (t)\left(  \frac{d}{dt} \right)^n \right)
  = \Res_0\,  dtf^{(n+1)} (t) g^{(m)} (t) \, .
\end{equation}
The principal gradation of $\D$ lifts to $\hat{\D}$ by letting
$\deg C=0$.  Note also that the restriction of the cocycle $\Psi$ 
to $\D^-$ is zero.

For each $s \in \CC$ one defines a Lie algebra homomorphism
$\varphi_s :\D \to \tilde{g\ell}_{\infty}$ (via the action of
$\D$ on $t^s \CC [t,t^{-1}]$) by
\begin{equation}
  \label{eq:3}
  \varphi_s (t^k f (D)) = \sum_{j \in \ZZ}
  f (-j+s) E_{j-k,j} \, .
\end{equation}
This homomorphism lifts to a homomorphism of central extension
$\hat{\varphi}_s : \hat{\D} \to \widehat{g\ell}_{\infty}$ by 
\begin{eqnarray}
    \label{eq:4}
\hat{\varphi}_s (t^k f (D))
   &=& \varphi_s (t^k f(D)) \hbox{ if } k \neq 0 \, , \\
\nonumber
\hat{\varphi}_s (e^xD) 
   &=& \varphi_s (e^{xD}) - \frac{e^{sx}-1}{e^x-1} \, ,
   \hat{\varphi}_s (C) =C \, .
\end{eqnarray}
More generally, for each $m \in \ZZ_+$ one defines a homomorphism 
$\varphi^{[m]}_s : \D \to \tilde{g\ell}_{\infty}^{[m]}$ by
\begin{equation}
  \label{eq:5}
  \varphi^{[m]}_s (t^k f (D)) = \sum_{j \in \ZZ}
  f (-j+s+u) E_{j-k,j} = \sum_{i=0}^m \sum_{j \in \ZZ}
  \frac{f^{(i)} (-j+s)}{i\!} u^i E_{j-k,j}  \, ,
\end{equation}
which lifts to $\hat{\varphi}_s^{[m]} : \hat{\D} \to
\widehat{g\ell}^{[m]}_{\infty}$ in a similar way.  One of the main
results of \cite{KR1} is the following.

\begin{lemma}
  \label{lem:3}
For each $i=1, \ldots ,r$, pick a collection $m_i \in \ZZ_+$,
$s_i \in \CC$, $\vec{\lambda}_i \in (\CC^{\infty})^{m_i+1}$,
$\vec{c}_i \in \CC^{m_i+1}$, such that $s_i-s_j \notin \ZZ$ for $
i \neq j$.  Then the $\oplus^r_{i=1} \widehat{g\ell}^{[m_i]}_{\infty}$-module
$\otimes^r_{i=1} L^{[m_i]} (\vec{\lambda}_i,\vec{c}_i)$ remains
irreducible when restricted to $\hat{\D}$ via the embedding
$\oplus^r_{i=1} \hat{\varphi}^{[m_i]}_{s_i} : \hat{\D} \to
\oplus^r_{i=1} \widehat{g\ell}^{[m_i]}_{\infty}$.  All irreducible
quasifinite highest weight $\hat{\D}$-modules are obtained in
this way.
\end{lemma}

\begin{proof}[Proof of Theorem~\ref{th:2}]
  Note that for $j \geq 1$ one has:
  \begin{equation}
    \label{eq:6}
    \D^-_j = \{ t^{-j} f (D) | f (0) = f(1) = \ldots = f(j-1)=0 
    \} \, .
  \end{equation}
Hence $\D^-_j$ has finite codimension in $\D_j$ and therefore the 
quasifiniteness of a $\D^-$-module $L(\Delta ;\D^-)$
implies the quasifiniteness of any of the $\hat{\D}$-modules 
$L(\Delta,c;\hat{\D})$.  Due to 
Lemma~\ref{lem:3}
, $L (\Delta ,c ;\hat{\D})$ is a tensor product of the
$\widehat{g\ell}^{[m]}_{\infty}$-modules $L^{[m]}(\vec{\lambda},
\vec{c})$ on which $\hat{\D}$ acts via the embedding
$\hat{\varphi}^{[m]}_s$ defined by (\ref{eq:3}) and (\ref{eq:4}).

It is clear from Theorem~\ref{th:1} that all non-trivial modules
$L^{[m]} (\vec{\lambda}_i, \vec{c}_i)$ have infinite growth (by
choosing an appropriate subalgebra isomorphic to
$g\ell_{+\infty}$ in $g\ell_{\infty}$).

Recall that for any quasifinite $\hat{\D}$-module one can extend
the action of $\hat{\D}_j$ for $j \neq 0$ to $\hat{\D}^\O_j$,
where $\O$ is the algebra of all holomorphic functions on $\CC$
\cite{KR1}, in other words, in (\ref{eq:3}) and (\ref{eq:4}) one can
take any $f \in \O$ if $j \neq 0$.  The same holds for $\D^-$,
except that for $j \geq 1$, $f$ must obey conditions in
(\ref{eq:6}).  We apply this to the $\hat{\D}$-module
$L(\widehat{g\ell}^{[m]}_{\infty}; \vec{\lambda},\vec{c})$ on which $\hat{\D}$ 
acts via $\hat{\varphi}^{[m]}_s$.

Choosing $f \in \O$ such that for all $j \in \ZZ$:
\begin{displaymath}
  f (-j+s)=\delta_{rj}, \,\,
  f^{(i)} (-j+s)=0 \hbox{ if } i=1,\ldots,m \, ,
\end{displaymath}
we see from (\ref{eq:3}) that all operators $E_{r+1,r}$ lie in
the image of $\hat{\varphi}^{[m]}_s (\D^{-\O})$, except for
$E_{1,0}$ when $s=0$ (here we use (\ref{eq:6}) for $j=1$).
Hence, when restricted to $\D^-$, the module $L^{[m]}
(\vec{\lambda},\vec{c})$ remains irreducible, provided that $s
\neq 0$.  Thus, if $L(\Delta ; \D^-)$ has finite growth, then $L
(\Delta ; \hat{\D})=L^{[m]}(\vec{\lambda},\vec{c})$ on which
$\hat{\D}$ acts via the embedding $\hat{\varphi}^{[m]}_0$.

Choosing $f \in \O$ to vanish in all $j \in \ZZ$ up to $m$\st{th} 
derivative except for $i$\st{th} derivative $(0<i \leq m)$ at
$j=-r$, we see that all operators $u^i E_{r+1,r}$ with $0<i \leq
m$ lie in the image of $\hat{\varphi}^{[m]}_s (\D^{-\O})$.

Suppose that the $m$\st{th} coordinate of $\vec{\lambda}_r$ is
non-zero, and that $m>0$.  Then $v:= (u^m E_{r+1,r})^N v_{\vec{\lambda}} 
\neq 0$ for all $N >0$.  But
\begin{displaymath}
E_{rr}v = (-N+ \lambda^0_r)v, \,
E_{r+1,r+1}v=(N+\lambda^0_{r+1})v.
\end{displaymath}
Therefore, restricting to
the subalgebra of $g\ell_{\infty}$ consisting of matrices
$(a_{ij})_{i,j \leq r}$ or $(a_{ij})_{i,j \geq r+1}$ we conclude
by Theorem~\ref{th:1}, that $L(\widehat{g\ell}^{[m]}_{\infty}; 
\vec{\lambda},\vec{c})$ is
either trivial or is of infinite growth.

Thus, the only possibility that remains is $s=m=0$.  As has been
already shown, the image of $\hat{\varphi}_s (\D^{-\O})$ contains 
all $E_{r+1,r}$ except for $E_{1,0}$, hence it contains all
operators from $g\ell_{-\infty} \oplus g\ell_{+\infty}$.
Therefore, by Theorem~\ref{th:1}
, the highest weight of a finite growth $\D^-$-module must be the
same as one of the $\D^-$-modules $L^+(\lambda^+) \otimes L^-
(\lambda^-)$ with $\lambda^{\pm} \in \Par^{\pm}$.

\end{proof}

Given  two partitions $\lambda^{\pm} \in \Par^{\pm}$, we
denote by $L (\lambda^+,\lambda^-)$ the $\D^-$-module, obtained
by restriction via $\varphi_0$ from the
$\tilde{g\ell}_{+\infty}\oplus \tilde{g\ell}_{-\infty}$-module\break
$L^+(\lambda^+) \otimes L^- (\lambda^-)$.  Now we shall
construct the $\D^-$-modules $L(\lambda^+,\lambda^-)$ explicitly.

Consider the $\D^-$-module $\CC [t,t^{-1}]$.  Then $\CC [t]$ is
its maximal submodule (which is irreducible).  Hence the
$\D^-$-module 
\begin{equation}\label{eq:V}
	V:=\CC [t,t^{-1}]/\CC [t]
\end{equation}
is irreducible.  It is 
clear that this is the highest weight $\D^-$-module of growth $1$ 
with a highest weight vector $t^{-1} + \CC [t]$.  It is immediate 
to deduce that $V$ is isomorphic to $L(\omega_1,0)$ where
$\omega_1=(1,0,\ldots) \in \Par^+$.

Likewise, the $\D^-$-module $\CC [t]^*=\oplus_{j \in \ZZ_+} (\CC
t^j)^*$ is an irreducible highest weight module of growth~$1$
with a highest weight vector $1^*$, hence it is isomorphic to
$L(0,\omega_{-1})$, where $\omega_{-1}=(\ldots ,0,-1) \in \Par^-$.  We
denote this $\D^-$-module by $V'$.

As in the Schur-Weyl theory, the $\D^-$-module $T^M (V) \otimes
T^N (V')$ has a natural decomposition as
$(\D^-,S_M \times S_N)$-modules:
\begin{displaymath}
  T^M (V) \otimes T^N (V')
  = \oplus_{\substack{\lambda^{\pm} \in \Par^{\pm}\\
    |\lambda^+|=M\\
    |\lambda^-|=N}} (V_{\lambda^+} \otimes V'_{\lambda^-})
  \otimes (U_{\lambda^+} \otimes U_{\lambda^-})
\end{displaymath}
where $U_{\lambda^+ \,(\hbox{resp. } \lambda^-)}$ denotes the irreducible 
$S_{M\, (\hbox{resp. }N)}$-module corresponding to the partition
$\lambda^+$ (resp. $\lambda^-$).

\begin{lemma}
  \label{lem:4}
The $\D^-$-modules $V_{\lambda^+} \otimes V'_{\lambda^-}$ are irreducible.
\end{lemma}

\begin{proof}
  As in the proof of Theorem~\ref{th:2}, we extend the action of
  $\D^-$ on $V_{\lambda^+} \otimes V'_{\lambda^-}$ to
  $\D^{-\O}_j$ for each $j \neq 0$, to obtain that any
  $\D^-$-submodule of $V_{\lambda^+} \otimes V'_{\lambda^-}$ is a 
  submodule over $g\ell_{+\infty} \oplus g\ell_{-\infty}$.  But, by 
  Schur--Weyl theory, the $g\ell_{+\infty} \oplus
  g\ell_{-\infty}$-module $V_{\lambda^+} \otimes V'_{\lambda^-}$
  is irreducible, which completes the proof.
\end{proof}

Thus, we have proved

\begin{theorem}
  \label{th:3}
  The $\D^-$-module $L (\lambda^+ ,\lambda^-)$ is
  isomorphic to $V_{\lambda^+} \otimes V'_{\lambda^-}$ for any
  pair $\lambda^{\pm} \in \Par^{\pm}$.

\end{theorem}

\begin{remark*}
  Considering $\lambda = (\lambda^-,\lambda^+) \in \CC^{\infty}$
  we may say that irreducible highest weight $\D^-$-modules of
  finite growth are parameterized by non-increasing sequences of integers
  $(\lambda_j)_{j \in \ZZ} \in \CC^{\infty}$ with the exception
  that $\lambda_0 \leq \lambda_1$.  Equivalently, letting $m_i
  =\lambda_i-\lambda_{i+1}$ we may say that these modules are
  parameterized by sequences of non-negative integers $(m_i)_{i
    \in \ZZ \backslash \{0\}}$, all but finite numbers of which
  are zero.
\end{remark*}

Recall that the extended annihilation algebra $\Lie^- (gc_1)$ for 
$gc_1$ is isomorphic to the direct sum of the Lie algebra $\D^-$
and the $1$-dimensional Lie algebra $\displaystyle{\CC (\partial
  + \tfrac{d}{dt})}$ and that conformal modules for a Lie
conformal algebra coincide with the conformal modules over the
associated extended annihilation algebra \cite{K2}.

Given a module $M$ over a Lie conformal algebra $R$ and $\alpha \in 
\CC$, we may construct the $\alpha$-twisted module $M_{\alpha}$
by replacing $\partial$ by $\partial +\alpha$ in the formulas for 
action of $R$ on $M$.  Theorems~\ref{th:2} and~\ref{th:3} and the 
above remarks imply

\begin{theorem}
  \label{th:4}
  The $gc_1$-modules $L(\lambda^+, \lambda^-)_{\alpha}$, where
  $\lambda^{\pm} \in \Par^{\pm}$, $\alpha \in \CC$, exhaust all
  irreducible conformal $gc_1$-modules of finite growth.

\end{theorem}

\begin{corollary*}
  The $gc_1$-modules $\CC [\partial]_{\alpha}$ and $\CC
  [\partial]^*_{\alpha}$, where $\alpha \in \CC$, exhaust all
  finite irreducible $gc_1$-modules.
  
\end{corollary*}

\begin{remark*}
  It is straightforward to generalize Theorems~\ref{th:2}
  and~\ref{th:3} to the case of $N \times N$ matrix differential
  operators and hence Theorem~\ref{th:4} to the case of $gc_N$.
  In particular the $gc_N$-modules $\CC [\partial]^N_{\alpha}$
  and $(\CC [\partial]^{N*})_{\alpha}$
  , where $\alpha \in \CC$, exhaust all finite irreducible
  $gc_N$-modules. This is an result of Kac, Radul and Wakimoto.
Moreover, these authors completely described all finite $gc_N$-modules,
which amonted to prove a complete reducibility result for finite
modules over the anninhilation algebra  (see \cite{K4}).
It is an open question whether a similar complete reducibility
result holds in the case of finite growth modules
\end{remark*}

\vskip 1cm


\section{Irreducible finite growth $gc_{1,\, x}\,$-modules}\label{sec:gc_1,x}

\vskip 1cm

The results of this section are almost the same as the previous one, as well as 
the proofs. Therefore, we will skip the details.

Let $\D_0$ (resp. $\D_0^-$) be the Lie subalgebra of $\D $ (resp. $\D^-$) of all 
regular differential operators on $\CC^*$ (resp. $\CC$) that kill constants. 
That is, $\D_0$ consists of linear combinations of  elements of the form $t^k D 
f(D)$, where $f$ is a polynomial. Denote by $\hat{\D_0}$ the corresponding 
central extension. These algebras inherits the $\ZZ$-gradation from $\hat\D$.

In this section, we will  need the representation theory of the Lie algebra
$\D^-_0$.  

Given a sequence of complex numbers $\Delta = (\Delta_1 ,
\Delta_2 , \ldots)$ we define the highest weight module $L(\Delta 
; \D^-_0)$ over $\D^-_0$ as the (unique) irreducible module that has
a non-zero vector $v_{\Delta}$ with the following properties:
\begin{displaymath}
  (\D^-_0)_j v_{\Delta} =0 \hbox{ for } j<0 , \,\,
    D^n v_{\Delta} = \Delta_n v_{\Delta}
    \hbox{ for }n \in \NN \, .
\end{displaymath}
The principal gradation of $\D^-_0$ induces the principal gradation $L (\Delta ; 
\D^-_0)$.

Quasifinite modules over $\D^-_0$ can be constructed as follows. The 
$\D^-_0$-modules $\CC
[t,t^{-1}]/\CC [t]$ and $\CC [t]/\CC$,  give us an embedding of 
$\D^-_0$ in $\tilde{g\ell}_{+\infty}$ and
$\tilde{g\ell}_{-\infty}$, respectively, hence an embedding of
$\D^-_0$ in $\tilde{g}\ell_{+\infty} \oplus
\tilde{g}\ell_{-\infty}$.  All these embeddings respect the
principal gradations.  Now take $\lambda^{\pm} \in \CC^{\pm
  \infty}$ and consider the $\tilde{g\ell}_{+\infty} \oplus
\tilde{g\ell}_{-\infty}$-module $L^+ (\lambda^+) \otimes L^-
(\lambda^-)$.  

 The same argument as in \cite{KR1}, gives us the
following.

\begin{lemma}
  \label{lem:8.1}
When restricted to $\D^-_0$, the module $L^+ (\lambda^+) \otimes
L^- (\lambda^-)$ remains irreducible.

\end{lemma}

It follows immediately that $L^+ (\lambda^+) \otimes L^-
(\lambda^-)$ is an irreducible highest weight module over $\D^-_0$, 
which is obviously quasifinite.  

We have the following theorem.

\begin{theorem}
  \label{th:8.2}
  The $\D^-_0$-modules $L^+ (\lambda^+) \otimes L^- (\lambda^-)$,
  where $\lambda^{\pm} \in \Par^{\pm}$, 
  exhaust all quasifinite irreducible highest weight
  $\D^-_0$-modules that have finite growth.
\end{theorem}

The proof of Theorem~\ref{th:8.2}  is the same as Theorem~\ref{th:2}, but in 
this case we reduce the problem to the  representation theory of the universal 
central extension 
$\hat{\D_0}$ of $\D_0$ that was   developed in 
\cite{AFMO} and \cite{KL}.

Let $s\in \ZZ $ and denote by  $\widehat{gl}^{[m]}_{\infty,s}$  the Lie 
subalgebra of $\widehat{gl}^{[m]}_{\infty}$ generated by $C$ and $\{u^l E_{ij} | 
0\leq l\leq m,  i\neq s \hbox{ and }  j\neq s \}$. Observe that 
$\widehat{gl}^{[m]}_{\infty,s}$ is naturally isomorphic to 
$\widehat{gl}^{[m]}_{\infty}$. Let $p_s:\widehat{gl}^{[m]}_{\infty}\to 
\widehat{gl}^{[m]}_{\infty,s}\to\widehat{gl}^{[m]}_{\infty}$ be the projection 
map composed with this isomorphism. If $s\notin \ZZ$, we also denote by 
$\hat{\varphi}^{[m]}_s$ the homomorphism  (\ref{eq:5}) restricted to 
$\hat{\D_0}$. If $s\in \ZZ$, we  redefine $\widehat{\varphi}^{[m]}_s$ by the 
homomorphism $p_s\circ \widehat{\varphi}^{[m]}_s:\hat\D_0 \to 
\widehat{gl}^{[m]}_{\infty,s}$. 

 In this case, we should replace Lemma \ref{lem:3} by one of the 
results of \cite{KL} (see also \cite{AFMO}):

\begin{lemma}
  \label{lem:8.3}
For each $i=1, \ldots ,r$, pick a collection $m_i \in \ZZ_+$,
$s_i \in \CC$, $\vec{\lambda}_i \in (\CC^{\infty})^{m_i+1}$,
$\vec{c}_i \in \CC^{m_i+1}$, such that $s_i-s_j \notin \ZZ$ for $
i \neq j$.  Then the $\oplus^r_{i=1} \widehat{g\ell}^{[m_i]}_{\infty}$-module
$\otimes^r_{i=1} L^{[m_i]} (\vec{\lambda}_i,\vec{c}_i)$ remains
irreducible when restricted to $\hat{\D_0}$ via the embedding
$\oplus^r_{i=1} \hat{\varphi}^{[m_i]}_{s_i} : \hat{\D_0} \to
\oplus^r_{i=1} \widehat{g\ell}^{[m_i]}_{\infty}$.  All irreducible
quasifinite highest weight $\hat{\D_0}$-modules are obtained in
this way.
\end{lemma}

\begin{proof}[Proof of Theorem~\ref{th:8.2}] The proof is the same as 
Theorem~\ref{th:2} but have to use Lemma~\ref{lem:8.3},  and in the case $s=0$ 
one should use the redefined $\hat{\varphi}^{[m]}_0$.

\end{proof}

Given  two partitions $\lambda^{\pm} \in \Par^{\pm}$, the $\D^-$-module  $L 
(\lambda^+,\lambda^-)$ that is obtained
by restriction via $\varphi_0$ from the
$\tilde{g\ell}_{+\infty}\oplus \tilde{g\ell}_{-\infty}$-module\break
$L^+(\lambda^+) \otimes L^- (\lambda^-)$, remains irreducible as a 
$\D_0^-$-module.  The construction of  the $\D^-_0$-modules 
$L(\lambda^+,\lambda^-)$ is the same as before, and Lemma \ref{lem:4} and 
Theorem \ref{th:3} holds for $\D^-_0$.

In this case, the extended annihilation algebra $\Lie^- (gc_{1,x})$ for 
$gc_{1,x}$ is isomorphic to the direct sum of the Lie algebra $\D^-_0$
and the $1$-dimensional Lie algebra $\displaystyle{\CC (\partial
  + \tfrac{d}{dt})}$. Theorems~\ref{th:8.2} and~\ref{th:3} and the 
above remarks imply

\begin{theorem}
  \label{th:8.4}
  The $gc_{1,x}$-modules $L(\lambda^+, \lambda^-)_{\alpha}$, where
  $\lambda^{\pm} \in \Par^{\pm}$, $\alpha \in \CC$, exhaust all
  irreducible conformal $gc_{1,x}$-modules of finite growth.

\end{theorem}

\begin{corollary*}
  The $gc_{1,x}$-modules $\CC [\partial]_{\alpha}$ and $\CC
  [\partial]^*_{\alpha}$, where $\alpha \in \CC$, exhaust all
  finite irreducible $gc_{1,x}$-modules.
  
\end{corollary*}

\vskip 1cm


\section{Irreducible finite growth $oc_1$-modules}\label{sec:oc1}

\vskip 1cm

Now, consider the anti-involution $\sigma$ on $\D$ defined by (cf. \cite{KWY})
\begin{displaymath}
  \sigma (t) = t \, , \qquad \,\, \sigma \left( \frac{d}{dt}\right) = -  
\frac{d}{dt}\, .
\end{displaymath}
Denote by $\D_\sigma$ the fixed  subalgebra of $\D$ by $-\sigma$, namely: 
$\D_\sigma=\{a\in\D \, | \, \sigma(a)=-a\, \}$. This subalgebra corresponds to 
the Lie algebra denoted by $\D^+$ in \cite{KWY}. Let $\hat\D_\sigma= \D_\sigma + 
\CC C$ denote the central extension given by the restriction of the cocycle 
(\ref{eq:cocycle-d}) on $\D$.

We are interested in representation theory of the Lie subalgebra
$\D^-_\sigma=\D^-\cap \hat\D_\sigma$ of regular differential operators on $\CC$ 
that are invariant by $-\sigma$. Both subalgebras inherit a $\ZZ$-gradation from 
$\D$, since $\sigma$ preserve the principal $\ZZ$-gradation of $\D$, and we have 
$\D_\sigma= \oplus_{j\in \ZZ}(\D_\sigma)_j$, where
  \begin{equation}
    \label{eq:8}
    (\D_\sigma)_j  = \{ t^j g(D+(j+1)/2) \, |\, g(w)\in \CC[w]\hbox{ is odd 
}\,\} .
\end{equation}
In the case of $(\D^-_\sigma)_j$, we need to add condition (\ref{eq:6}) for 
$j<0$. 
 
Similarly, we have the corresponding subalgebras of $\D^\O$, denoted by 
$\D_\sigma^\O$ and $\D_\sigma^{- \O}$.

As in the case of $\D^-$, given a sequence of complex numbers $\Delta = 
\{\Delta_n\}_{n\in \NN_{odd}}$,
 we define the highest weight module $L(\Delta 
; \D^-_\sigma)$ over $\D^-_\sigma$ as the (unique) irreducible module that has
a non-zero vector $v_{\Delta}$ with the following properties:
\begin{displaymath}
  (\D^-_\sigma)_j v_{\Delta} =0 \hbox{ for } j<0 , \,\qquad
    (D+1/2)^n v_{\Delta} = \Delta_n v_{\Delta}
    \hbox{ for }n \in \NN_{odd} \, .
\end{displaymath}
The principal gradation of $\D^-_\sigma$ induces the principal gradation $L 
(\Delta ; \D^-_\sigma)
=\oplus_{j \in \ZZ_+} L_j$ such that $L_0 =\CC v_{\Delta}$.  The
module $L (\Delta ; \D^-_\sigma)$ is called \emph{quasifinite} if $\dim
L_j <\infty$ for all $j \in \ZZ_+$.

Quasifinite modules over $\D^-_\sigma$ can be constructed as follows.
The  $\D^-_\sigma$-module $\CC
[t,t^{-1}]/\CC [t]$   gives us an embedding of 
$\D^-_\sigma$ in $\tilde{g\ell}_{+\infty}$.  This embedding respect the
principal gradations.  Now take $\lambda^+  \in \CC^{+
  \infty}$ and consider the $\tilde{g\ell}_{+\infty}$-module $L^+ (\lambda^+ )$ 
introduced in (\ref{eq:1}).  The same argument as in \cite{KR1}, gives us the
following.

\begin{lemma}
  \label{lem:6}
When restricted to $\D^-_\sigma$, the module $L^+ (\lambda^+)$ remains 
irreducible.

\end{lemma}

Therefor $L^+ (\lambda^+)$ is an irreducible quasifinite highest weight module 
over $\D^-_\sigma$, 
and   it is easy to see that we have:
\begin{eqnarray*}
  \Delta_n &=& \sum_{j \geq 1} (-j+1/2)^n \lambda^+_j
   \, , \quad n \in\NN_{odd} \, ,\\
\noalign{\hbox{so that}}\\
\Delta (x): &=& \sum_{n \in\NN_{odd}} \Delta_n x^n /n! 
   = \sum_{j \geq 1}\lambda^+_j \, 2\, \sinh{((-j+1/2)x)}\, .
\end{eqnarray*}

We shall prove the following theorem.

\begin{theorem}
  \label{th:6}
  The $\D^-_\sigma$-modules $L^+ (\lambda^+)  $,
  where $\lambda^{+} \in \Par^{+}$, 
  exhaust all quasifinite irreducible highest weight
  $\D^-_\sigma$-modules that have finite growth.
\end{theorem}

The basic idea of the proof of Theorem~\ref{th:6} is the same as
in Theorem~\ref{th:2}:  to reduce the problem to the well developed (in
\cite{KWY}) representation theory of the universal central extension 
$\hat{\D}_\sigma$.

Recall that the homomorphism $\hat\varphi_s^{[m]} :\hat\D\to 
\widehat{g\ell}_\infty^{[m]}$ defined in (\ref{eq:5}) lift to a homomorphism  
$\hat\varphi_s^{[m]} :\hat\D^\O \to \widehat{g\ell}_\infty^{[m]}$. Now, the 
restriction  $\hat\varphi_s^{[m]} :\hat\D^\O_\sigma \to 
\widehat{g\ell}_\infty^{[m]}$ to $\hat\D^\O_\sigma$  is surjective iff $s\notin 
\ZZ/2$, and in the other cases, using (\ref{eq:8}), we have that (see \cite{KWY} 
for details)
\begin{equation}
  \label{eq:9}
  \hat\varphi_0^{[m]} : \hat\D_\sigma^\O \to {d}_\infty^{[m]} \, , \qquad \quad 
 \hat\varphi_{-1/2}^{[m]} : \hat\D_\sigma^\O \to {b}_\infty^{[m]}
\end{equation}
are surjective homomorphisms. Now, let us consider the restriction to 
$\D_\sigma^{- \O}$. Since the constrains given by (\ref{eq:6}) do not affect the 
case $s\neq 0$, we still have that 
$
\hat\varphi_s^{[m]} :\D_\sigma^{- \O} \to \widehat{g\ell}_\infty^{[m]}
$ ($s\notin \ZZ/2$)  and 
$
 \hat\varphi_{-1/2}^{[m]} : \D_\sigma^{- \O} \to {b}_\infty^{[m]}
$ are surjective.

  One of the main results of \cite{KWY} is the following.

\begin{lemma}
  \label{lem:7}
For each $i=1, \ldots ,r$, pick a collection $m_i \in \ZZ_+$,
$s_i \in \CC$, $\vec{\lambda}_i \in (\CC^{\infty})^{m_i+1}$,
$\vec{c}_i \in \CC^{m_i+1}$, such that  $s_i\in\ZZ$ implies $s_i=0$, 
$s_i\in\frac 1 2 +\ZZ$ implies $s_i=-\frac 1 2$, and $s_i-s_j \notin \ZZ$ for 
$i \neq j$.  Then the $\oplus^r_{i=1} \fg^{[m_i]}$-module
$\otimes^r_{i=1} L(\fg^{[m_i]};\vec{\lambda}_i,\vec{c}_i)$ remains
irreducible when restricted to $\hat{\D}_\sigma$ via the embedding
$\oplus^r_{i=1} \hat{\varphi}^{[m_i]}_{s_i} : \hat{\D}_\sigma \to
\oplus^r_{i=1} \fg^{[m_i]}$, where $\fg^{[m_i]}= \widehat{g\ell}_\infty^{[m_i]}$ 
(resp. ${b}_\infty^{[m_i]}$ or ${d}_\infty^{[m_i]}$) if $s_i\notin \ZZ/2$ (resp. 
$s_i=-\frac 1 2$ or $s_i=0$).  
All irreducible
quasifinite highest weight $\hat{\D}_\sigma$-modules are obtained in
this way.
\end{lemma}

\begin{proof}[Proof of Theorem~\ref{th:6}]
 The proof is similar to that of Theorem~\ref{th:2}.  Due to 
Lemma~\ref{lem:7}, Theorem~\ref{th:5} and (\ref{eq:9}), it is easy to see that  
if $L(\Delta ; \D^-_\sigma)$ has finite growth, then $L
(\Delta ; \hat{\D}_\sigma)=L(d_\infty^{[m]}; \vec{\lambda},\vec{c})$ on which
$\hat{\D}_\sigma$ acts via the embedding $\hat{\varphi}^{[m]}_0$.

Choosing $f \in \O_{odd}$ to vanish in all $j \in \ZZ$ up to $m$\st{th} 
derivative except for $i$\st{th} derivative $(0<i \leq m)$ at
$j=-r$, we see that all operators $u^i E_{r+1,r}- (-u)^i E_{-r+1,-r}$ with $0<i 
\leq
m$ lie in the image of $\hat{\varphi}^{[m]}_0 (\D^{-\O}_\sigma)$.

Suppose that the $m$\st{th} coordinate of $\vec{\lambda}_r$ is
non-zero, and that $m>0$.  Then $v:= (u^m E_{r+1,r}- (-u)^i E_{-r+1,-r})^N 
v_{\vec{\lambda}} 
\neq 0$ for all $N >0$.  But
\begin{displaymath}
E_{r+1,r+1}v=(N+\lambda^0_{r+1})v.
\end{displaymath}
As in Theorem~{\ref{th:2}}, restricting to
the subalgebra of $d_\infty^{[m]}$ isomorphic to $g\ell_{+\infty}$ consisting of 
matrices $(a_{i,j}- a_{1-j,1-i})_{i,j \geq r+1}$ we conclude
by Theorem~\ref{th:1}, that $L(d^{[m]}_\infty; \vec{\lambda},\vec{c})$ is
either trivial or is of infinite growth.

Thus, the only possibility that remains is $s=m=0$.  As has been
already shown, the image of $\hat{\varphi}_s (\D^{-\O}_\sigma)$ contains 
all $E_{r+1,r}-E_{1-r,-r}$ for all $r\neq 0$, hence it contains all
operators from $d_\infty \cap \left(g\ell_{-\infty} \oplus 
g\ell_{+\infty}\right)\simeq g\ell_{+\infty}$.
Therefore, by Theorem~\ref{th:1}, the highest weight of a finite growth 
$\D^-_\sigma$-module must be the
same as one of the $\D^-_\sigma$-modules $L^+(\lambda^+)$ with $\lambda^{+} \in 
\Par^{+}$.

\end{proof}

Now we shall
construct the $\D^-_\sigma$-modules $L^+(\lambda^+)$ explicitly. The
$\D^-$-module $V=\CC [t,t^{-1}]/\CC [t]$ defined in (\ref{eq:V}), viewed as a 
$\D^-_\sigma$-module, remains  irreducible.  This is the highest weight 
$\D^-_\sigma$-module of growth $1$ 
 isomorphic to $L^+(\omega_1)$ where
$\omega_1=(1,0,\ldots) \in \Par^+$.

Observe that the $\D^-_\sigma$-module $\CC [t]^*=\oplus_{j \in \ZZ_+} (\CC 
t^j)^*$ is isomorphic to
$L^+(\omega_{1})$.

As in the Schur-Weyl theory, the $\D^-_\sigma$-module $T^M (V)$ has a natural 
decomposition as
$(\D^-_\sigma,S_M)$-modules:
\begin{displaymath}
  T^M (V)
  = \oplus_{\substack{\lambda^{+} \in \Par^{+}\\
    |\lambda^+|=M}} V_{\lambda^+}
  \otimes U_{\lambda^+}
\end{displaymath}
where $U_{\lambda^+}$ denotes the irreducible 
$S_{M}$-module corresponding to the partition
$\lambda^+$.

\begin{lemma}
  \label{lem:8}
The $\D^-_\sigma$-modules $V_{\lambda^+}$ are irreducible.
\end{lemma}

\begin{proof}
  As in the proof of Theorem~\ref{th:6}, we extend the action of
  $\D^-_\sigma$ on $V_{\lambda^+}$ to
  $(\D^{-\O}_\sigma)_j$ for each $j \neq 0$, to obtain that any
  $\D^-_\sigma$-submodule of $V_{\lambda^+}$ is a 
  submodule over $g\ell_{+\infty}\, \left(\simeq d_\infty \cap (g\ell_{+\infty} 
\oplus g\ell_{-\infty})\right)$.  But, by 
  Schur--Weyl theory, the $g\ell_{+\infty}$-module $V_{\lambda^+}$
  is irreducible, which completes the proof.
\end{proof}

Thus, we have proved

\begin{theorem}
  \label{th:7}
The $\D^-_\sigma$-module $T^M (V)$ has the following decomposition  as 
\linebreak
$(\D^-_\sigma,S_M )$-modules
\begin{displaymath}
  T^M (V)
  = \oplus_{\substack{\lambda^{+} \in \Par^{+}\\
    |\lambda^+|=M }} 
    L^+({\lambda^+}) 
  \otimes U_{\lambda^+} 
\end{displaymath}
where $U_{\lambda^+ }$ denotes the irreducible 
$S_{M}$-module corresponding to the partition
$\lambda^+$.
\end{theorem}

\begin{remark*}
  Considering $\lambda^+ \in \CC^{+\infty}$
  we may say that irreducible highest weight $\D^-_\sigma$-modules of
  finite growth are parameterized by non-increasing sequences of integers
  $(\lambda_j)_{j \in \NN} \in \CC^{+\infty}$.  Equivalently, letting $m_i
  =\lambda_i-\lambda_{i+1}$ we may say that these modules are
  parameterized by sequences of non-negative integers $(m_i)_{i
    \in \NN }$, all but finite numbers of which
  are zero.
\end{remark*}

Recall that the extended annihilation algebra $\Lie^- (oc_1)$ for 
$oc_1$ is isomorphic to the direct sum of the Lie algebra $\D^-_\sigma$
and the $1$-dimensional Lie algebra $\displaystyle{\CC (\partial
  + \tfrac{d}{dt})}$ and that conformal modules for a Lie
conformal algebra coincide with the conformal modules over the
associated extended annihilation algebra \cite{K2}.

  Theorems~\ref{th:6} and the 
above remarks imply

\begin{theorem}
  \label{th:8}
  The $oc_1$-modules $L^+(\lambda^+)_{\alpha}$, where
  $\lambda^{+} \in \Par^{+}$, $\alpha \in \CC$, exhaust all
  irreducible conformal $oc_1$-modules of finite growth.

\end{theorem}

\begin{corollary*}
  The $gc_1$-modules $L^+(\lambda^+)$, where
  $\lambda^{+} \in \Par^{+}$, remain irreducible when restricted
to $oc_1$.
  
\end{corollary*}

\begin{corollary*}
  The $oc_1$-modules $\CC [\partial]_{\alpha}$, where $\alpha \in \CC$, exhaust 
all
  finite irreducible $oc_1$-modules.
  
\end{corollary*}


\vskip 1cm


\section{Irreducible finite growth $spc_1$-modules}

\vskip 1cm

Now, consider the anti-involution $\bar\sigma$ on $\D_0$ defined by 
\begin{displaymath}
  \bar\sigma (t^k D f(D)) = - t^k Df(-D-k).
\end{displaymath}
This anti-involution was studied  by Bloch \cite{B} in connection with the 
values of $\zeta$-function.

Denote by $\D_{0,\bar\sigma}$ the Lie subalgebra of $\D_0$ fixed by 
$-\bar\sigma$.  Let $\hat\D_{0,\bar\sigma}= \D_{0,\bar\sigma} + \CC C$ denote 
the central extension given by the restriction of the cocycle on $\D$.

We are interested in representation theory of the Lie subalgebra
$\D^-_{0,\bar\sigma}=\D^-\cap \hat\D_{0,\bar\sigma}$ of regular differential 
operators on $\CC$ that kills constants and are invariant by $-\bar\sigma$. Both 
subalgebras inherit a $\ZZ$-gradation from $\D_0$, since $\bar\sigma$ preserve 
the principal $\ZZ$-gradation of $\D_0$: $\D_{0,\bar\sigma}= \oplus_{j\in 
\ZZ}(\D_{0,\bar\sigma})_j$, where
  \begin{equation}
    \label{eq:4.3}
    (\D_{0,\bar\sigma})_j  = \{ t^j D\, g(D+\frac{j}{2}) \, |\, g(w)\in 
\CC[w]\hbox{ is even }\,\} .
\end{equation}
In the case of $(\D^-_{0,\bar\sigma})_j$, we need to add condition (\ref{eq:6}) 
for $j<0$. 
 
Similarly, we have the corresponding subalgebras of $\D^\O$, denoted by 
$\D_{0,\bar\sigma}^\O$ and $\D_{0,\bar\sigma}^{- \O}$.

As in the case of $\D^-$, given a sequence of complex numbers $\Delta = 
\{\Delta_n\}_{n\in \NN_{odd}}$,
 we define the highest weight module $L(\Delta 
; \D^-_{0,\bar\sigma})$ over $\D^-_{0,\bar\sigma}$ as the (unique) irreducible 
module that has
a non-zero vector $v_{\Delta}$ with the following properties:
\begin{displaymath}
  (\D^-_{0,\bar\sigma})_j v_{\Delta} =0 \hbox{ for } j<0 , \,\qquad
    D^n v_{\Delta} = \Delta_n v_{\Delta}
    \hbox{ for }n \in \NN_{odd} \, .
\end{displaymath}
The principal gradation of $\D^-_{0,\bar\sigma}$ induces the principal gradation 
$L (\Delta ; \D^-_{0,\bar\sigma})
=\oplus_{j \in \ZZ_+} L_j$ such that $L_0 =\CC v_{\Delta}$.  The
module $L (\Delta ; \D^-_{0,\bar\sigma})$ is called \emph{quasifinite} if $\dim
L_j <\infty$ for all $j \in \ZZ_+$.

As in the previous section, the  $\D^-_{0,\bar\sigma}$-module $\CC
[t,t^{-1}]/\CC [t]$   gives us an embedding of 
$\D^-_{0,\bar\sigma}$ in $\tilde{g\ell}_{+\infty}$.  This embedding respect the
principal gradations.  Now take $\lambda^+  \in \CC^{+
  \infty}$ and consider the $\tilde{g\ell}_{+\infty}$-module $L^+ (\lambda^+ )$ 
introduced in (\ref{eq:1}).  The same argument as in \cite{KR1}, gives us the
following.

\begin{lemma}
  \label{lem:4.3}
When restricted to $\D^-_{0,\bar\sigma}$, the quasifinite module $L^+ 
(\lambda^+)$ remains irreducible.

\end{lemma}

We shall prove the following theorem.

\begin{theorem}
  \label{th:4.4}
  The $\D^-_{0,\bar\sigma}$-modules $L^+ (\lambda^+)  $,
  where $\lambda^{+} \in \Par^{+}$, 
  exhaust all quasifinite irreducible highest weight
  $\D^-_{0,\bar\sigma}$-modules that have finite growth.
\end{theorem}

The basic idea of the proof of Theorem~\ref{th:4.4} is the same as
in Theorem~\ref{th:2}:  to reduce the problem to the recently developed (in
\cite{BL}) representation theory of the universal central extension 
$\hat{\D}_{0,\bar\sigma}$.

Recall that the homomorphism $\hat\varphi_s^{[m]} :\hat\D\to 
\widehat{g\ell}_\infty^{[m]}$ defined in (\ref{eq:5}) lift to a homomorphism  
$\hat\varphi_s^{[m]} :\hat\D^\O \to \widehat{g\ell}_\infty^{[m]}$. Now, the 
restriction  $\hat\varphi_s^{[m]} :\hat\D^\O_{0,\bar\sigma} \to 
\widehat{g\ell}_\infty^{[m]}$ to $\hat\D^\O_{0,\bar\sigma}$  is surjective iff 
$s\notin \ZZ/2$, and in the other case, using (\ref{eq:4.3}), we have that (see 
[BL] for details)
\begin{equation}
  \label{eq:4.5}
  \hat\varphi_s^{[m]} : \hat\D_{0,\bar\sigma}^\O \to {c}_\infty^{[m]} \, , 
\qquad \quad s\in \ZZ/2
\end{equation}
is a  surjective homomorphism. Now, let us consider the restriction to 
$\D_{0,\bar\sigma}^{- \O}$. Since the constrains given by (\ref{eq:6}) do not 
affect the case $s\neq 0$, we still have that 
$
\hat\varphi_s^{[m]} :\D_{0,\bar\sigma}^{- \O} \to \widehat{g\ell}_\infty^{[m]}
$ ($s\notin \ZZ/2$)  and 
$
 \hat\varphi_{-1/2}^{[m]} : \D_{0,\bar\sigma}^{- \O} \to {c}_\infty^{[m]}
$ are surjective.

  One of the main results of \cite{BL} is the following.

\begin{lemma}
  \label{lem:4.6}
For each $i=1, \ldots ,r$, pick a collection $m_i \in \ZZ_+$,
$s_i \in \CC$, $\vec{\lambda}_i \in (\CC^{\infty})^{m_i+1}$,
$\vec{c}_i \in \CC^{m_i+1}$, such that  $s_i\in\ZZ$ implies $s_i=0$, 
$s_i\in\frac 1 2 +\ZZ$ implies $s_i=-\frac 1 2$, and $s_i-s_j \notin \ZZ$ for 
$i \neq j$.  Then the $\oplus^r_{i=1} \fg^{[m_i]}$-module
$\otimes^r_{i=1} L(\fg^{[m_i]};\vec{\lambda}_i,\vec{c}_i)$ remains
irreducible when restricted to $\hat{\D}_{0,\bar\sigma}$ via the embedding
$\oplus^r_{i=1} \hat{\varphi}^{[m_i]}_{s_i} : \hat{\D}_{0,\bar\sigma} \to
\oplus^r_{i=1} \fg^{[m_i]}$, where $\fg^{[m_i]}= \widehat{g\ell}_\infty^{[m_i]}$ 
(resp. ${c}_\infty^{[m_i]}$) if $s_i\notin \ZZ/2$ (resp. $s_i=-\frac 1 2$ or 
$s_i=0$).  
All irreducible
quasifinite highest weight $\hat{\D}_{0,\bar\sigma}$-modules are obtained in
this way.
\end{lemma}

\begin{proof}[Proof of Theorem~\ref{th:4.4}]
 The proof is similar to that of Theorem~\ref{th:2}.  Due to 
Lemma~\ref{lem:4.6}, Theorem~\ref{th:4.2} and (\ref{eq:4.5}), it is easy to see 
that  if $L(\Delta ; \D^-_{0,\bar\sigma})$ has finite growth, then $L
(\Delta ; \hat{\D}_{0,\bar\sigma})=L(c_\infty^{[m]}; \vec{\lambda},\vec{c})$ on 
which
$\hat{\D}_{0,\bar\sigma}$ acts via the embedding $\hat{\varphi}^{[m]}_0$.

Choosing $f \in \O_{odd}$ to vanish in all $j \in \ZZ$ up to $m$\st{th} 
derivative except for $i$\st{th} derivative $(0<i \leq m)$ at
$j=-r$, we see that all operators $u^i E_{r+1,r} + (-u)^i E_{-r+1,-r}$ with $0<i 
\leq
m$ lie in the image of $\hat{\varphi}^{[m]}_0 (\D^{-\O}_{0,\bar\sigma})$.

Suppose that the $m$\st{th} coordinate of $\vec{\lambda}_r$ is
non-zero, and that $m>0$.  Then $v:= (u^m E_{r+1,r} + (-u)^i E_{-r+1,-r})^N 
v_{\vec{\lambda}} 
\neq 0$ for all $N >0$.  But
\begin{displaymath}
E_{r+1,r+1}v=(N+\lambda^0_{r+1})v.
\end{displaymath}
As in Theorem~{\ref{th:2}}, restricting to
the subalgebra of $c_\infty^{[m]}$ isomorphic to $g\ell_{+\infty}$ consisting of 
matrices $(a_{i,j}- (-1)^{i+j}a_{1-j,1-i})_{i,j \geq r+1}$ we conclude
by Theorem~\ref{th:1}, that $L(c^{[m]}_\infty; \vec{\lambda},\vec{c})$ is
either trivial or is of infinite growth.

Thus, the only possibility that remains is $s=m=0$.  As has been
already shown, the image of $\hat{\varphi}_s (\D^{-\O}_{0,\bar\sigma})$ contains 
all $E_{r+1,r} + E_{1-r,-r}$ for all $r\neq 0$, hence it contains all
operators from $c_\infty \cap \left(g\ell_{-\infty} \oplus 
g\ell_{+\infty}\right)\simeq g\ell_{+\infty}$.
Therefore, by Theorem~\ref{th:1}, the highest weight of a finite growth 
$\D^-_{0,\bar\sigma}$-module must be the
same as one of the $\D^-_{0,\bar\sigma}$-modules $L^+(\lambda^+)$ with 
$\lambda^{+} \in \Par^{+}$.

\end{proof}

As in the previous section, we can 
construct the $\D^-_{0,\bar\sigma}$-modules $L^+(\lambda^+)$ explicitly. The
$\D^-$-module $V=\CC [t,t^{-1}]/\CC [t]$ defined in (\ref{eq:V}), viewed as a 
$\D^-_{0,\bar\sigma}$-module, remains  irreducible.  This is the highest weight 
$\D^-_{0,\bar\sigma}$-module of growth $1$ 
 isomorphic to $L^+(\omega_1)$ where
$\omega_1=(1,0,\ldots) \in \Par^+$.

Observe that the $\D^-_{0,\bar\sigma}$-module $\CC [t]^*=\oplus_{j \in \ZZ_+} 
(\CC t^j)^*$ is isomorphic to
$L^+(\omega_{1})$.

As in the Schur-Weyl theory, the $\D^-_{0,\bar\sigma}$-module $T^M (V)$ has a 
natural decomposition as
$(\D^-_{0,\bar\sigma},S_M)$-modules:
\begin{displaymath}
  T^M (V)
  = \oplus_{\substack{\lambda^{+} \in \Par^{+}\\
    |\lambda^+|=M}} V_{\lambda^+}
  \otimes U_{\lambda^+}
\end{displaymath}
where $U_{\lambda^+}$ denotes the irreducible 
$S_{M}$-module corresponding to the partition
$\lambda^+$.

\begin{lemma}
  \label{lem:4.7}
The $\D^-_{0,\bar\sigma}$-modules $V_{\lambda^+}$ are irreducible.
\end{lemma}

\begin{proof}
  As in the proof of Theorem~\ref{th:6}, we extend the action of
  $\D^-_{0,\bar\sigma}$ on $V_{\lambda^+}$ to
  $(\D^{-\O}_{0,\bar\sigma})_j$ for each $j \neq 0$, to obtain that any
  $\D^-_{0,\bar\sigma}$-submodule of $V_{\lambda^+}$ is a 
  submodule over $g\ell_{+\infty}\, \left(\simeq c_\infty \cap (g\ell_{+\infty} 
\oplus g\ell_{-\infty})\right)$.  But, by 
  Schur--Weyl theory, the $g\ell_{+\infty}$-module $V_{\lambda^+}$
  is irreducible, which completes the proof.
\end{proof}

Thus, we have proved

\begin{theorem}
  \label{th:4.8}
The $\D^-_{0,\bar\sigma}$-module $T^M (V)$ has the following decomposition  as
$(\D^-_{0,\bar\sigma},S_M )$-modules
\begin{displaymath}
  T^M (V)
  = \oplus_{\substack{\lambda^{+} \in \Par^{+}\\
    |\lambda^+|=M }} 
    L^+({\lambda^+}) 
  \otimes U_{\lambda^+} 
\end{displaymath}
where $U_{\lambda^+ }$ denotes the irreducible 
$S_{M}$-module corresponding to the partition
$\lambda^+$.
\end{theorem}

Recall that the extended annihilation algebra $\Lie^- (spc_1)$ for 
$spc_1$ is isomorphic to the direct sum of the Lie algebra $\D^-_{0,\bar\sigma}$
and the $1$-dimensional Lie algebra $\displaystyle{\CC (\partial
  + \tfrac{d}{dt})}$ and that conformal modules for a Lie
conformal algebra coincide with the conformal modules over the
associated extended annihilation algebra \cite{K2}.

  Theorems~\ref{th:4.4} and the 
above remarks imply

\begin{theorem}
  \label{th:4.9}
  The $spc_1$-modules $L^+(\lambda^+)_{\alpha}$, where
  $\lambda^{+} \in \Par^{+}$, $\alpha \in \CC$, exhaust all
  irreducible conformal $spc_1$-modules of finite growth.

\end{theorem}

\begin{corollary*}
  The $gc_1$-modules $L^+(\lambda^+)$, where
  $\lambda^{+} \in \Par^{+}$, remain irreducible when restricted to $spc_1$.
  
\end{corollary*}

\begin{corollary*}
  The $spc_1$-modules $\CC [\partial]_{\alpha}$, where $\alpha \in \CC$, exhaust 
all
  finite irreducible $spc_1$-modules.
  
\end{corollary*}

\vskip 1cm

\noindent{\bf Acknowledgment.} C. Boyallian and J. Liberati were supported in
 part by
Conicet, ANPCyT, Agencia Cba Ciencia,  
Secyt-UNC and Fomec (Argentina). Special thanks go to MSRI (Berkeley) for the 
hospitality during our stay there, where part of this work was done.




\begin{thebibliography}{11}


\bibitem[1]{AFMO} H. Awata, M. Fukuma,  Y. Matsuo and  S. Odake, {\em 
Subalgebras
of $W\sb {1+\infty}$ and their quasifinite
representations}, J. Phys. A {\bf 28} (1995), 105-112.


\bibitem[2]{B} S. Bloch, {\em Zeta values and differential operators on
the circle}, J. Algebra {\bf 182} (1996),   476-500.  

\bibitem [3]{BKL} C. Boyallian, V. G. Kac and J. I.  Liberati, {\em 
On the classification of subalgebras of Cend$_N$ and $gc_N$},  Journal of 
Algebra (2002).
math-ph/0203022
 
\bibitem [4]{BL} C. Boyallian and J. I.  Liberati, {\em 
Representations of a symplectic type subalgebra of $W_\infty$},  submitted   
(2002).


\bibitem[5]{DeK}  A. De Sole and V.~G.~Kac, {\em  Subalgebras of $gc_N$ and 
Jacobi polynomials}, 
 Canadian Math. Bull. (2002). math-ph/0112028.



\bibitem[6]{K1}
V.~G.~Kac,
{\textit{Infinite-dimensional Lie algebras}},
3rd edition, Cambridge University Press, Cambridge, 1990.

\bibitem[7]{K2}
V.~G.~Kac,
{\textit{Vertex algebras for beginners}}, 
University Lecture Series, 10. American Mathematical Society, Providence, RI,
1996. Second edition 1998.



\bibitem[8]{K4}
V.~G.~Kac,
{\textit{Formal distribution algebras and conformal algebras}},
in Proc.
XIIth International Congress of
Mathematical Physics (ICMP '97) (Brisbane), 80--97, 
Internat. Press, Cambridge, MA, 1999;
{\texttt{q-alg/9709027}}.



\bibitem [9]{KL}  V. G. Kac and  J. I. Liberati, {\em 
Unitary quasifinite representations of $W_\infty$}, Letters in Math. Phys. {\bf 
53} (2000), 11-27.

\bibitem[10]{KR1}  V.~G.~Kac  and A. Radul, {\em  Quasifinite highest weight
modules  over the Lie algebra of differential operators on the circle}, 
 Comm. Math. Phys. {\bf 157} (1993), 429-457.



\bibitem[11]{KWY}  V.~G.~Kac, W. Wang and C. Yan,  {\em Quasifinite
representations of classical Lie subalgebras of  $W_{1+\infty}$}
Adv. Math. {\bf 139} (1998),  56-140.


\end{thebibliography}
\end{document}